# On the Zeros of the Complex Fourier Transforms of a Class of Exponential Functions

Jeremy C Williams, London, England
E-mail: universaltutors@aol.com

Abstract: A class theorem is presented and proved: the complex Fourier transforms of a certain class of exponential functions have all their zeros on the real line. A class of basis functions is first considered, and the class is then extended by the method of convolutions.

In this paper we are going to consider the complex Fourier transforms of the exponential kernel $e^{-q(t)}$, where the following conditions apply:

(i) q(t) is an even function of t which is real for all real values of t

(ii) q(it) is either a polynomial which has all its roots on the real line, or the uniform limit on compact sets in the complex number field of such a polynomial.

It is well-known that (ii) corresponds to q(it) belonging to the Laguerre-Pólya class of functions (denote by L-P).

If there are no non-zero roots $\beta_n$ of q(it) = 0, then q(it) = c. $t^N.e^{-\mu t^2+\gamma t}$ (1.0)

and if there are non-zero roots $\beta_n$ of q(it) = 0, then q(it) = c. $t^N.e^{-\mu t^2+\gamma t}.\prod_n (1-t^2/\beta_n^2).e^{-t/\beta_n}$ (1.1)

We note $\gamma, \beta_n \in \Re$, the $\beta_n$ are the non-zero roots of q(it) = 0, $\mu \geq 0$, N is a non-negative integer, and that the product in (1.1) may or may not be infinite.

Consider the case when there are non-zero roots $\beta_n$. Replacing t with t/i in equation (1.1) shows that:

q(t) = c/$i^N.t^N$. $e^{(\mu t^2)+t(\gamma-\Sigma 1/\beta_n)/i}$ $.\prod_n (1+t^2/\beta_n^2)$ (1.2)

As q is an even function of t, it follows that for every non-zero root $\beta_n$ of q(t) = 0, there is a root of opposite sign -$\beta_n$, hence $\sum_n (1/\beta_n)$ = 0.

Hence equation (1.2) reduces to:

q(t)=c/$i^N$ . $t^N.e^{\gamma t/i}$. $e^{\mu t^2}$. $\prod_n (1+t^2/\beta_n^2)$ =c/$i^N$ . $e^{\mu t^2} t^N.(\cos(\gamma t)-i.\sin(\gamma t))\,.\prod_n (1+t^2/\beta_n^2)$ (1.3)

We note that the above expression can be split into the difference of 2 parts:

q(t) = ( $c/i^N.e^{\mu t^2}.t^N.\cos(\gamma t).\prod_n (1+t^2/\beta_n^2)$ )- ( $c/i^N.e^{\mu t^2}.t^N.\sin(\gamma t).\prod_n (1+t^2/\beta_n^2)$ )

If N is odd then the first part is odd and the second part is even for all real values of t. As q(t) is an even function in t, the odd part must be zero, forcing cos( $\gamma t$ ) = 0 for all real t, which is impossible, so we have obtained a contradiction. Therefore N must be even, so the first part is even and the second odd. So sin( $\gamma t$) = 0 for all real values of t which means we must have $\gamma$ = 0.

Denote N by 2m, and let m > 0 for simplicity because then q(0) = 0 . We can now see that c must be real, because of (i), and that:

q(t) = $(-1)^m.c.t^{2m}.e^{\mu t^2}.\prod_n (1+t^2/\beta_n^2)$ , and writing $(-1)^m.c$ as k, we obtain:

$$q(t) = k.\, t^{2m}.e^{\mu t^2}.\prod_n (1+t^2/\beta_n^2) \qquad (1.4)$$

In the case when there are no roots, the working above is shortened and we just obtain:

$$q(t) = k.\, t^{2m}.e^{\mu . t^2} \qquad (1.5)$$

Therefore, let us work with equation (1.4), and take the product to be 1 if there are no non-zero roots $\beta_n$ . We have shown that if we require q(t) to be an even function of t which is real for all real values of t, and q(it) to be either a polynomial, or the uniform limit on compact sets in the complex number field of such a polynomial, such that q(it) has all its roots on the real line, that this leads to q(t) being expressible in the form of equation (1.4).

Let us also require that k > 0, and that there are positive constants T, $\alpha$ such that for real t, q(t) > $t^{2+\alpha}$ for t ≥ T (this condition is equivalent to there being a non-zero term in $t^4$ or higher in the series expansion of q(t) about t = 0. For example, if the coefficient of $t^4$ in this series expansion is $a_4$ > 0, then we see q(t) ≥ $a_4.t^4$ , and so q(t) > $t^{2+\alpha}$ for t > $(1/a_4)^{1/(2-\alpha)}$ , 0 < $\alpha$ < 2).

Define $G(z) = \int_{-\infty}^{\infty} e^{-q(t)}.e^{izt}\,dt$ = R($\sigma$,w) + i I($\sigma$,w), where z = w - i$\sigma$, and $\sigma$,w, R, I $\in \Re$ (2.0)

**Theorem 1** G(z) as defined above has only real zeros.

**Proof of Theorem 1**

The proof proceeds by defining a sequence of functions $g_n(t)$ which converge to g(t) = $e^{-q(t)}$, and their complex Fourier transforms $G_n(z)$. It is proved that $G_n(z)$ has only real zeros and converges absolutely and uniformly to G(z) on compact subsets in the complex number field.

The convolution $h_n(t) = g_n(t) \circ e^{-t^4}$ and its complex Fourier transform $H_n(z)$ are then considered. It is proved that $H_n(z)$ is a function of order less than 2 in z, with an infinite number of real zeros, and consequently that $H_n(z)$ is a Laguerre - Pólya function. It is then shown why, as a result of the above, all the zeros of G(z) are real:

The preceding definition of q(t) ensures that it is strictly increasing above 0 for t > 0. Hence q(t) is a 1:1 function on $(0,\infty)$, and for any natural number n, there is exactly one real root of the equation q(t) = n in the interval $(0,\infty)$. Call this root $\lambda_n$ (by the evenness of q, we see q(-$\lambda_n$) = n as well ). We also note that as q(t) $\to\infty$ as t $\to\infty$, it follows that $\lambda_n \to \infty$ as n $\to\infty$ (each $\lambda_n$ is finite).

Define $g_n(t) = \chi_{[-\lambda_n,\lambda_n]}(t)\,.\,y_n(t)$ where $y_n(t) = (\frac{n}{n+1})\,.(1-q(t)/n)^{n+1}$ for all real t (2.1)

where $\chi_{[a,b]}(t) = 1 \quad a \le t \le b$

$\quad\quad 0 \quad$ otherwise

Now also define a sequence of functions $G_n(z) = \int_{-\infty}^{\infty} g_n(t).e^{izt}.dt$ (2.2)

In 1927 (reference I), Professor Pólya proved that if F(t) is defined for real values of t and equals $\overline{F(-t)}$, and there exist positive real constants A, $\alpha$ such that $|F(t)|$ < A $e^{-|t|^{(2+\alpha)}}$ for all real values of t, then if all the zeros of $\int_{-\infty}^{\infty} F(t).e^{izt}dt$ are real, so are the zeros of $\int_{-\infty}^{\infty} \phi(t).F(t).e^{izt}.dt$, provided that $\phi(it)$ is either a polynomial with only real roots, or the uniform limit on compact sets of such. Call this result (R1).

We see that $\int_{-\infty}^{\infty} \chi_{[-\lambda_n,\lambda_n]}(t).e^{izt}dt$. = $\int_{-\lambda_n}^{\lambda_n} e^{izt}.dt$ = $[e^{izt}/(iz)]_{t=-\lambda_n}^{t=\lambda_n}$ = $(e^{iz\lambda_n} - e^{-iz\lambda_n})$/(iz)

which equals zero if and only if $e^{2iz\lambda_n} = e^{i.2k\pi}$, for any integer k, in other words only when z is real and equal to k$\pi/\lambda_n$. Now set $\chi_{[-\lambda_n,\lambda_n]}(t)$ = F(t) in (R1) above, we see that $\chi_{[-\lambda_n,\lambda_n]}(t)$ has compact support, so it easily obeys the condition $|F(t)|$ < A.$e^{-|t|^{(2+\alpha)}}$ for some positive constants A, $\alpha$. Recalling that L-P is closed under differentiation, we observe that as q(it) lies in L-P, then so does $q'(it)$, and so $q'(it)$ is either a polynomial with only real roots, or the uniform limit on compact sets of such, as in (R1) above.

Hence by (R1), we see that $\int_{-\infty}^{\infty} q'(t).\chi_{[-\lambda_n,\lambda_n]}(t).e^{izt}dt = \int_{-\lambda_n}^{\lambda_n} q'(t).e^{izt}dt$ only has real zeros

Now consider $\int_{-\infty}^{\infty} (1-q(t)/n)^1.\chi_{[-\lambda_n,\lambda_n]}(t).e^{izt}dt = \int_{-\lambda_n}^{\lambda_n} (1-q(t)/n)^1.e^{izt}.dt$

Integration by parts shows that this is equal to $[(1-q(t)/n)^1.e^{izt}/iz]_{-\lambda_n}^{\lambda_n}$ - $\int_{-\lambda_n}^{\lambda_n} -\frac{q'(t)}{n}.\frac{e^{izt}}{iz}.dt$

= 0 + $\frac{1}{izn}$. $\int_{-\lambda_n}^{\lambda_n} q'(t).e^{izt}dt$ and by the preceding result has only real zeros.

Now set F(t) = $(1-q(t)/n)^1.\chi_{[-\lambda_n,\lambda_n]}(t)$. By (R1) we see that $\int_{-\infty}^{\infty} q'(t).(1-q(t)/n)^1.\chi_{[-\lambda_n,\lambda_n]}(t).e^{izt}dt$

= $\int_{-\lambda_n}^{\lambda_n} q'(t).(1-q(t)/n)^1.\, e^{izt}dt$ has only real zeros.

Consider $\int_{-\infty}^{\infty} (1-q(t)/n)^2.\chi_{[-\lambda_n,\lambda_n]}(t).e^{izt}.dt$ = $\int_{-\lambda_n}^{\lambda_n} (1-q(t)/n)^2.e^{izt}.dt$.

Again, integration by parts shows us a useful result, namely that:

$\int_{-\lambda_n}^{\lambda_n} (1-q(t)/n)^2.e^{izt}.dt = [(1-q(t)/n)^2.e^{izt}/iz]_{-\lambda_n}^{\lambda_n}$ - $\frac{-2}{izn}.\int_{-\lambda_n}^{\lambda_n} q'(t).(1-q(t)/n)^1.e^{izt}.dt$

= 0 + $\frac{2}{izn}\int_{-\lambda_n}^{\lambda_n} q'(t).(1-q(t)/n)^1.e^{izt}.dt$ which by the above also has only real zeros.

Continuing this argument shows that $\int_{-\infty}^{\infty} (1-q(t)/n)^{(n+1)}.\chi_{[-\lambda_n,\lambda_n]}.e^{izt}.dt$ and hence $G_n(z)$ $= \; (\frac{n}{n+1}).\int_{-\infty}^{\infty}(1-q(t)/n)^{(n+1)}.\chi_{[-\lambda_n,\lambda_n]}.e^{izt}.dt \; = \; \int_{-\lambda_n}^{\lambda_n}(\frac{n}{n+1}).(1-q(t)/n)^{(n+1)}.e^{izt}.dt$ has only real roots.

Call this result (R2)

We will now see that $G_n(z)$ converges absolutely and uniformly to G(z):

Consider the disc $|z| \le$ M, and choose $\varepsilon > 0$, where $\varepsilon$ is arbitrarily small.

We see that $|G(z)-G_n(z)| = \left|\int_{-\infty}^{\infty} e^{-q(t)}.e^{izt}dt - (\frac{n}{n+1})\int_{-\infty}^{\infty}\chi_{[-\lambda_n,\lambda_n]}(t).(1-q(t)/n)^{(n+1)}.e^{izt}dt\right|$

$$\le \int_{-\infty}^{\infty}\left|(e^{-q(t)}-g_n(t)).e^{izt}\right|.dt \; \le \int_{-\infty}^{\infty}\left|e^{-q(t)}-g_n(t)\right|e^{|z|t}.dt \; \le \int_{-\infty}^{\infty}\left|e^{-q(t)}-g_n(t)\right|.e^{Mt}.dt$$

$$\le 2\int_{0}^{\infty}\left|e^{-q(t)}-g_n(t)\right|.e^{Mt}.dt$$

therefore $|G(z)-G_n(z)| \; \le 2(\int_{0}^{t_1}\left|e^{-q(t)}-g_n(t)\right|.e^{Mt}.dt + \int_{t_1}^{\infty}\left|e^{-q(t)}-g_n(t)\right|.e^{Mt}.dt$ ) for some $t_1 > 0$ (2.4)

Call $2\int_{0}^{t_1}\left|e^{-q(t)}-g_n(t)\right|.e^{Mt}.dt = J_n$ (2.5)

and $2\int_{t_1}^{\infty}\left|e^{-q(t)}-g_n(t)\right|.e^{Mt}.dt \; = \; K_n$ (2.6)

Regarding $K_n$, we begin by noting that when $|x| < 1$, the series expansion about x = 0 of ln(1-x)

is $-(x+\frac{x^2}{2}+\frac{x^3}{3}+\frac{x^4}{4}+.....) \; \le \; -x$, and replacing x with $\frac{q(t)}{n}$ yields the relation:

$\ln(1-q(t)/n) \le$ - q(t)/n for $|q(t)/n| < 1$, so we see that n $\ln(1-q(t)/n) \le -q(t)$

and hence that $(1-q(t)/n)^n \le e^{-q(t)}$ for $t \in [0,\lambda_n)$

Considering the case t = $\lambda_n$ shows that this inequality holds on the closed interval $[0,\lambda_n]$

ie that $(1-q(t)/n)^n \le e^{-q(t)}$ for $t \in [0,\lambda_n]$ (2.7)

We also observe that $0 \le (1-q(t)/n) \le 1$ for $t \in [0,\lambda_n]$ and so $(1-q(t)/n)^n \ge 0$. We also note that $0 < \frac{n}{n+1} < 1$ for n > 1.

Multiplying the left hand side of equation (2.4) by (1-q(t)/n) and n/(n+1), we see that:

$0 \le \frac{n}{n+1}.(1-q(t)/n)^{n+1} \le e^{-q(t)}$ for $t \in [0,\lambda_n]$ (2.8)

We can now see that if t $\in [0,\lambda_n]$ then both $e^{-q(t)} - g_n(t) \geq 0$ and $g_n(t) \geq 0$, and hence that

$\left|e^{-q(t)} - g_n(t)\right|$ = $e^{-q(t)} - g_n(t) \leq e^{-q(t)}$, whereas if t > $\lambda_n$, then $g_n(t) = 0$, so $\left|e^{-q(t)} - g_n(t)\right|$

= $\left|e^{-q(t)}\right|$ = $e^{-q(t)}$, so we can see that from equation (2.6) it follows that :

$$K_n \leq 2\int_{t_1}^{\infty} e^{-q(t)}.e^{Mt}.dt \qquad (2.9)$$

As q(t) increases more quickly than $t^{2+\alpha}$ for $t \geq T$, the integral in (2.9) exists for any fixed finite value of M, and we can make the integral less than $\varepsilon/2$ by choosing $t_1$ to be sufficiently large. So $t_1$ is now chosen.

Regarding $J_n$, and recall that $y_n(t) = e^{-q(t)} - (\frac{n}{n+1}).(1-q(t)/n)^{n+1}$. We have just seen that $y_n(t) \geq 0$ for t $\in [0,\lambda_n]$, and we note that $y_n'(t) = -q'(t).(e^{-q(t)} - (1-q(t)/n)^n)$. From equation (1.5) we see that $q'(t)$ is zero for t = 0 and strictly positive for t > 0. By equation (2.7), $e^{-q(t)} - (1-q(t)/n)^n \geq 0$ and so $y_n'(t) \leq 0$ for t $\in [0,\lambda_n]$.

Now pick n > q($t_1$), so $\lambda_n > t_1$. As $y_n'(t) \leq 0$ and $y_n(t) \geq 0$ for t$\in [0,t_1]$, it follows that the maximum value of $y_n(t)$ on the interval $[0,t_1]$ is $y_n(0)$ = 1/(n+1).

Therefore equation (2.5) leads us to:

$$J_n \leq 2\int_0^{t_1} (\frac{1}{n+1}).e^{Mt}.dt = 2[\frac{e^{Mt}}{M(n+1)}]_0^{t_1} = 2\,(e^{Mt_1}\text{-}1)/(M(n+1)) < 2\,e^{Mt_1}/(Mn)$$

Choosing any n > Max[q($t_1$), $4e^{Mt_1}/(M\varepsilon)$] ensures that the above conditions are satisfied and that $J_n < \in/2$ and by (2.4) that $\left|G(z) - G_n(z)\right| < \varepsilon/2 + \varepsilon/2 = \varepsilon$ for all z such that $|z| \leq M$, and as $\varepsilon$ is arbitrarily small, we have proved the uniform and absolute convergence of $G_n(z)$ to G(z).

……………………………………………………..……………………………………………..

Now let us consider the convolutions $h_n(t) = g_n(t) \circ e^{-t^4}$ and h(t) = $e^{-q(t)} \circ e^{-t^4}$ defined by :

$$h_n(t) = \int_{-\infty}^{\infty} g_n(v).e^{-(t-v)^4}.dv \qquad (3.0)$$

and h(t) = $\int_{-\infty}^{\infty} e^{-q(t)}.e^{-(t-v)^4}.dv$ (3.1)

We see that $H_n(z)$, the complex Fourier transform of $h_n(t)$ defined by $H_n(z) = \int_{-\infty}^{\infty} h_n(t).e^{izt}.dt$ is then equal to $\iint_{\Re^2} (g_n(v).e^{-(t-v)^4}).e^{izt}.dvdt$ and after substituting v = $x$ and t = $x + y$ and applying Fubini’s Theorem, we see that:

$$H_n(z) = (\int_{-\infty}^{\infty} g_n(x).e^{izx}.dx).(\int_{-\infty}^{\infty} e^{-y^4}.e^{izy}.dy) = G_n(z).F_4(z) \text{ where } F_4(z) = \int_{-\infty}^{\infty} e^{-y^4}.e^{izy}.dy \quad (3.1)$$

Of course equation (3.1) is just the standard result that the Fourier transform of the convolution is the product of the individual Fourier transforms.

We saw earlier (R2) that $G_n(z)$ has only real zeros, and in 1923 (ref II) Pólya proved that for k $\geq 2$, all the zeros of $F_{2k}(z) = \int_{-\infty}^{\infty} e^{-y^{2k}}.e^{izy}.dy$ lie on the real line Im(z) = 0, and that there are an infinite number of them.

Hence $F_4(z)$, and consequently $H_n(z)$, have only real zeros, and further more an infinity of them. Call this result (R3).

**The order of $H_n(z)$:**

In the same paper previously mentioned (I), Pólya also proved that the <u>order</u> of the complex Fourier transform of f(t), where f(t) is a real function of t such that f(t) = $\overline{f(-t)}$ and there exist real positive constants A, $\alpha$ such that $|f(t)|$ < A.$e^{-|t|^{2+\alpha}}$ for all real values of t, is at most (2+$\alpha$)/(1+$\alpha$) < 2 . Call this result (R4).

Therefore, let us set f(t) = $g_n(t) = \chi_{[-\lambda_n,\lambda_n]}(t).(\frac{n}{n+1}).(1-q(t)/n)^{n+1}$. Outside the interval $[-\lambda_n,\lambda_n]$ we observe $g_n(t) = 0$, and inside, by equation (2.8), $0 \leq g_n(t) \leq e^{-q(t)}$.

As q(t) > $t^{2+\alpha}$ for t $\geq$ T, it follows that $e^{-q(t)} < e^{-t^{2+\alpha}}$ for t $\geq$ T, and by setting A = 2 Sup[$e^{-q(t)+t^{2+\alpha}}$] on the interval $0 \leq t \leq T$ , we see $e^{-q(t)} <$ A.$e^{-|t|^{2+\alpha}}$ and so by (R4) above we conclude the order of $G_n(z)$ is at most (2+$\alpha$ )/(1+$\alpha$ ) < 2

Incidentally, setting f(t) = $e^{-q(t)}$ instead of $g_n(t)$ shows us, by the same reasoning as that in the above paragraph, that the order of G(z) is at most (2+$\alpha$ )/(1+$\alpha$ ) < 2 as well.

Setting f(t) = $e^{-t^4}$ shows that the order of $F_4(z)$ is at most 4/3. (Pólya mentions in (II) that the order of $F_{2k}(z)$ is 2k/(2k-1), which we obtain from (R4) by writing $\alpha$ = 2k-2 ).

Therefore, there must exist constants B, $c_1, c_2$ > 0 such that :

$$|H_n(z)| = |G_n(z).F_4(z)| \leq Be^{c_1|z|^{\frac{2+\alpha}{1+\alpha}}}.e^{c_2|z|^{4/3}} \leq Be^{(c_1+c_2)|z|^{\rho}} \text{ ,where } \rho = Max[\frac{2+\alpha}{1+\alpha}, \frac{4}{3}], \quad \text{for all z}$$

We then see that the order of $H_n(z)$ is less than 2.

Note: another way of evaluating an upper limit for the order of $H_n(z)$ is to observe that $h_n(t) = \int_{-\lambda_n}^{\lambda_n} y_n(t).e^{-(t-v)^4}.dv$ and $|h_n(t)| \leq \int_{-\lambda_n}^{\lambda_n} |y_n(t)|.e^{-(t-v)^4}.dv$ and use the result that

$|y_n(t)| = \left|(\frac{n}{n+1}).(1-q(t)/n)^{n+1}\right| \le 1$ for t $\in [-\lambda_n, \lambda_n]$ so we must have

$|h_n(t)| \le \int_{-\lambda_n}^{\lambda_n} e^{-(v-t)^4}.dv = \int_{t-\lambda_n}^{t+\lambda_n} e^{-v^4}.dv$ and if we consider t > $2\lambda_n$ then $|h_n(t)| \le 2\,\lambda_n$ .Max[$e^{-v^4}$]

on the interval $[-\lambda_n, \lambda_n]$ ie $|h_n(t)| \le 2\lambda_n.e^{-(t-\lambda_n)^4} \le 2\lambda_n.e^{-(t/2)^4}$, so we see for large t that $h_n(t)$ falls off more quickly than $e^{-t^{2+\alpha}}$, and so satisfies the conditions in (R4).

We can then set f(t) = $h_n(t)$ and the result that the order of $H_n(z)$ is less than 2 follows.

## $H_n(z)$ and $G_n(z)$ are Laguerre- Pólya functions

As we saw in (R3) earlier, $H_n(z)$ has only real zeros, and an infinity of them, and we have just seen in the preceding section that the order of $H_n(z)$ is less than 2.

By Hadamard's factorisation theorem (see for example III), $H_n(z)$ has the product representation:

$$H_n(z) = \text{c.}\, z^N.e^{\gamma z}.\prod_{k=1}^{\infty}(1+z/w_k)e^{-z/w_k} = \text{c.}\, z^n.e^{(\gamma-\Sigma 1/w_j)z}.\prod_{k=1}^{\infty}(1+z/w_k) \qquad (3.2)$$

We will now see that this expression can be reduced to a rather simple form, using the same reasoning as in a previous paper by the same author (IV).

From the definitions of $G_n(z)$ and $F_4(z)$ in (2.2) and (3.1), we see that they are both even functions, and therefore so is $H_n(z) = G_n(z).F_4(z)$. Consequently, for every non-zero root $w_k$ of $H_n(z) = 0$, there is another root of opposite sign - $w_k$. So we see that $\sum_j 1/w_j = 0$.

Equation (3.2) then reduces to:

$$H_n(z) = c.z^N.e^{\gamma z}.\prod_{r=1}^{\infty}(1 - z^2/\alpha_r^2) \qquad (3.3)$$

where $w_1 = \alpha_1$, $w_2 = -\alpha_1$, $w_3 = \alpha_2$, $w_4 = -\alpha_2$ etc and the $\alpha_r$ are the positive roots of $H_n(z) = 0$ .

As $H_n(0) = G_n(0).F_4(0) = (\int_{-\lambda_n}^{\lambda_n} \frac{n}{n+1}.(1-g(t)/n)^{n+1}.dt).(\int_{-\infty}^{\infty} e^{-t^4}.dt)$ and 1-g(t)/n > 0 when $|t| < \lambda_n$,

we see $H_n(0) > 0$ and so N = 0 and c > 0. Therefore (3.3) reduces to:

$$H_n(z) = c.e^{\gamma z}.\prod_{r=1}^{\infty}(1 - z^2/\alpha_r^2) \qquad (3.4)$$

and as $H_n(z)$ is even in z, we see $\gamma = 0$ and so that:

$$H_n(z) = c.\prod_{r=1}^{\infty}(1 - z^2/\alpha_r^2), \text{ where } c > 0 \qquad (3.5)$$

Recalling that the order of $F_4(z)$ is at most 4/3, and that $F_4(z)$ has only real zeros, and furthermore an infinity of them, the same reasoning leads to the conclusion that:

$$F_4(z) = c_0.\prod_{k=1}^{\infty}(1 - z^2/\beta_k^2) \quad \text{where } c_0 = \int_{-\infty}^{\infty} e^{-t^4}.dt > 0 \qquad (3.6)$$

As $H_n(z) = G_n(z).F_4(z)$, *every zero of* $F_4(z)$ *is a zero of* $H_n(z)$ and so it follows from (3.5) and (3.6) that:

$$G_n(z) = \frac{H_n(z)}{F_4(z)} = \frac{c.\prod_{r=1}^{\infty}(1 - z^2/\alpha_r^2)}{c_0.\prod_{k=1}^{\infty}(1 - z^2/\beta_k^2)} = \frac{c}{c_0}.\prod_j (1 - z^2/(\alpha_j^*)^2) \qquad (3.7)$$

where the set $\{\alpha_j^*\}_{j\geq 1}$ is the set $\{\alpha_r\}_{r=1}^{r=\infty}$ with the set $\{\beta_k\}_{k=1}^{k=\infty}$ removed ie $\{\alpha_r\}_{r=1}^{r=\infty} \setminus \{\beta_k\}_{k=1}^{k=\infty}$.

…………………………………………………………………………………………………

In the same paper by this author previously mentioned (IV), we saw that for

$F_{2n}(z) = \int_{-\infty}^{\infty} e^{-t^{2n}}.e^{izt}dt$, where z = w - i$\sigma$,

$$|F_{2n}(z)|^2 = \frac{1}{2}\sum_{m\geq 0}\frac{\sigma^{2m}}{(2m)!}.A_{2m,2n}(w) \qquad (3.8a)$$

where $$A_{2m,2n}(w) = \frac{1}{i^{2m}}.\frac{\partial^{2m}}{\partial u^{2m}}\left[\iint_{\Re^2} e^{-(((t+X)/2)^{2n} + ((t-X)/2)^{2n}) + i(wX+ut)}.dtdX\right]_{u=0} \qquad (3.9a)$$

and so that $$A_{2m,2n}(w) = (-1)^m.2.\frac{\partial^{2m}}{\partial u^{2m}}\left[F_{2n}(u+w).F_{2n}(u-w)\right]_{u=0} \qquad (4.0a)$$

and we also saw that $$T_{N,m}(w) = (-1)^m.2.\frac{\partial^{2m}}{\partial u^{2m}}\left[P_N(u+w).P_N(u-w)\right]_{u=0} \geq 0 \qquad (4.1a)$$

where $P_N(z) =$ c. $\prod_{r=1}^{N}(1 - z^2/\alpha_r^2)$ and $F_{2n}(z) = c.\prod_{r=1}^{\infty}(1 - z^2/\alpha_r^2)$ (4.2a)

and by letting N $\to \infty$, that $A_{2m,2n} \geq 0$.

By the same reasoning, when we consider $G_n(z) = \int_{-\infty}^{\infty} g_n(t).e^{izt}.dt$, we obtain:

$$|G_n(z)|^2 = \frac{1}{2}.\sum_{m\geq 0}\frac{\sigma^{2m}}{(2m)!}.B_{2m,n}(w) \qquad (3.8b)$$

where $$B_{2m,n}(w) = \frac{1}{i^{2m}}.\frac{\partial^{2m}}{\partial u^{2m}}\left[\iint_{\Re^2} g_n((t+X)/2).g_n((t-X)/2).e^{i(wX+ut)}.dtdX\right]_{u=0} \qquad (3.9b)$$

and so $B_{2m,n}(w) = (-1)^m.2.\frac{\partial^{2m}}{\partial u^{2m}}\left[B_n(u+w).B_n(u-w)\right]_{u=0}$ (4.0b)

with $T_{N,m}(w) = (-1)^m.2.\frac{\partial^{2m}}{\partial u^{2m}}\left[P_N(u+w).P_N(u-w)\right]_{u=0} \geq 0$ (4.1b)

and $P_N(z) = \frac{c}{c_0}.\prod_{r=1}^{N}(1 - z^2/(\alpha_r^*)^2)$ and $G_n(z) = \frac{c}{c_0}.\prod_{r\geq 1}(1 - z^2/(\alpha_r^*)^2)$ (4.2b)

Whether the product representation for $G_n(z)$ in (4.2b) above is infinite or not, we see that $B_{2m,n}(w) \geq 0$.

We note that for any particular fixed values of n and w, it is not possible for $B_{2m,n}(w)$ to be zero for all values of m, otherwise by equation (3.8b), there would be a continuous line of zeros of $G_n(z)$ along the line w = constant, and this would lead to a contradiction as the zeros of holomorphic functions are isolated. Therefore, for any particular fixed values of n and w, at least one of the $B_{2m,n}(w)$ must be greater than zero.

Hence by equation (3.8b) we see that if we pick a fixed value of w on the w-axis and travel in the direction of increasing $\sigma$, then $|G_n(z)|^2$ *must increase monotonically.*

By the absolute uniform convergence of $G_n(z)$ to G(z), which ensures the limit of $G_n(z)$ is continuous, we immediately see that G(z) must have the same property. For this reason, it follows that all the zeros of G(z) lie on the line $\sigma = 0$.

In other words, all the zeros of G(z) are real.

QED.

…………………………………………………………………………………………………..

By way of example, if we set q(t) = a (cosh(t)-1), we see that $e^{-a}$ G(z) = $e^{-a}\int_{-\infty}^{\infty} e^{-a(\cosh(t)-1)}.e^{izt}.dt$

= $\int_{-\infty}^{\infty} e^{-a.\cosh(t)}.e^{izt}.dt$ and we note that q(t) = 2a $(\sinh(t/2))^2$ = $\frac{at^2}{2}.\prod_{m=1}^{\infty}(1+\frac{t^2}{4\pi^2 m^2})^2$ , which is of the form in equation (1.4). This case was covered in Pólya's 1927 paper (I), but it was proved directly, and not by a class rule.

…………………………………………………………………………………………………..

## Extending the class

We now turn our attention to the convolution of two basis functions $e^{-q_1(t)}$ and $e^{-q_2(t)}$ defined by:

$c_{[q_1,q_2]}(t) = \int_{-\infty}^{\infty} e^{-q_1(v)}.e^{-q_2(t-v)}.dv$ (4.3)

We see $c_{[q_1,q_2]}(-t)=\int_{-\infty}^{\infty} e^{-q_1(u)}.e^{-q_2(-t-u)}.du$ and as $q_2(t)$ is an even function, this equals

$\int_{-\infty}^{\infty} e^{-q_1(u)-q_2(t+u)}.du$ which on replacing u with –v becomes $\int_{-\infty}^{\infty} e^{-q_1(v)-q_2(t-v)}.dv$, in other words

$c_{[q_1,q_2]}(t)$. We have proved therefore, that the convolution is an even function of the variable t.

We note that $c_{[q_1,q_2]}(t)$ is always strictly positive as the integrand in (4.3) is strictly positive.

Therefore $c_{[q_1,q_2]}(t)$ must be expressible as $e^{-g_{1,2}(t)}$ where $g_{1,2}(t)$ is even in t.

Now, for convenience, consider t > 0 and differentiate the expression in (4.3) with respect to t. We obtain:

$c_{[q_1,q_2]}'(t)=\int_{-\infty}^{\infty} -q_2'(t-v).e^{-q_1(v)-q_2(t-v)}.dv$ which after the change of variable u = t – v becomes

$\int_{-\infty}^{\infty} -q_2'(u).e^{-q_1(t-u)-q_2(u)}.du$ and this expression may be broken up as:

$$(\int_0^{\infty}+\int_{-\infty}^{0})(-q_2'(u).e^{-q_1(t-u)-q_2(u)}).du=\int_0^{\infty} q_2'(u).e^{-q_2(u)}.(e^{-q_1(t+u)}-e^{-q_1(t-u)}).du$$

so $$c_{[q_1,q_2]}'(t)=\int_0^{t} q_2'(u).e^{-q_2(u)}.(e^{-q_1(t+u)}-e^{-q_1(t-u)}).du+\int_t^{\infty} q_2'(u).e^{-q_2(u)}.(e^{-q_1(t+u)}-e^{-q_1(u-t)}).du \quad (4.4)$$

For 0 <u <t  we see t + u > t - u  > 0 and so $e^{-q_1(t+u)}-e^{-q_1(t-u)}<0$ as $q_1(x)$ is strictly increasing on $(0,\infty)$. For u > t > 0 we see u + t > u – t > 0 and so $e^{-q_1(t+u)}-e^{-q_1(u-t)}<0$ as well.

It follows that both integrals in (4.4) are negative, and so we conclude that:

$c_{[q_1,q_2]}'(t)=-g_{1,2}'(t).e^{-g_{1,2}(t)}<0$ for all t > 0,

and so $g_{1,2}'(t)>0 \quad \forall t>0$.

We can see therefore that the plot of the convolution $c_{[q_1,q_2]}(t)$ against t is a symmetric bell-shaped curve which looks similar to that of $e^{-q(t)}$ where q(t) is as defined in (1.4).

By the standard result concerning the complex Fourier transform (CFT) of convolutions, the zeros of the CFT of $c_{[q_1,q_2]}(t)$ are exactly those of the individual CFT's of $e^{-q_1(t)}$ and $e^{-q_2(t)}$. Therefore all the zeros of the CFT's of the convolutions as defined in (4.3) are real.

*We therefore define our class of exponential kernels whose CFT's have only real zeros as the set* $e^{-q(t)}$ *where q(t) is as defined in (1.4), and their convolutions as defined by (4.3)*

…………………………………………………………………………………………………….

Author: Jeremy Williams, MA(Oxon)
E-mail: universaltutors@aol.com